\documentclass{article}

\usepackage{diagrams}
\usepackage{latexsym}
\def\text#1{\mbox{\rm #1\ }}
\def\ie{{\rm i.e.,\/}\ }

\def\cf{{\rm cf.\/}\ }
\def\id{\mbox{\it id\,}}
\def\one{\mbox{\rm 1}\hskip-2.8pt \mbox{\rm l}}
\def\ZZ{\mbox{\rm Z}\hskip-5pt \mbox{\rm Z}}
\def\RR{\mbox{\rm I}\hskip-2pt \mbox{\rm R}}
\def\CC{\mbox{\rm C}\hskip-5.5pt \mbox{l} \;}

\title{Finite dimensional quantum group covariant differential calculus
       on a complex matrix algebra
   \vspace{0.7cm} \\
}

\author{R. Coquereaux${}^1$\thanks{~Email: coque@cpt.univ-mrs.fr},
        A. O. Garc\'{\i}a${}^2$\thanks{~Email: ariel@cab.cnea.edu.ar},
        R. Trinchero${}^2$\thanks{~Email: trincher@cab.cnea.edu.ar} \\
\\
${}^1$ {\it Centre de Physique Th\'eorique - CNRS - Luminy, Case 907} \\
       {\it F-13288 Marseille Cedex 9 - France} \\
\\
${}^2$ {\it Instituto Balseiro and Centro At\'omico Bariloche} \\
       {\it CC 439 - 8400 - Bariloche - R\'{\i}o Negro - Argentina} \\
\\
}

\date{April 2, 1998}

\begin{document}

\thispagestyle{empty}
\begin{titlepage}

\maketitle

\vfill

\abstract{
Using the fact that the algebra $M_3(\CC)$ of $3 \times 3$ complex matrices
can be taken as a reduced quantum plane, we build a differential calculus
$\Omega(S)$ on the quantum space $S$ defined by the algebra
$C^\infty(M) \otimes M_3(\CC)$, where $M$ is a space-time manifold. This
calculus is covariant under the action and coaction of finite dimensional
dual quantum groups. We study the star structures on these quantum groups
and the compatible one in $M_3(\CC)$. This leads to an invariant scalar
product on the later space. We analyse the differential algebra
$\Omega(M_3(\CC))$ in terms of quantum group representations, and consider
in particular the space of $1$-forms on $S$ since its elements can be
considered as generalized gauge fields.
}

\vspace{1.2 cm}


\vspace{1.0cm}

\noindent Anonymous ftp or gopher: cpt.univ-mrs.fr

\vspace{0.5 cm}

\noindent {\tt math.QA/9804021}\\
\noindent CPT-98/P.3630 \\
\noindent IT-CNEA-CAB/2902798

\vspace*{0.3 cm}

\end{titlepage}


\section {Introduction}
\label{sec:introduction}

The formulation of physical theories in the framework of non commutative
geometry opens new possibilities and has produced very interesting
results. Physical models over a space $S$ described by the tensor product
of the commutative algebra of functions over a space-time manifold and the
non-commutative space whose ``algebra of functions'' is given by
${\mathcal M} = M_N(\CC)$, the algebra of $N \times N$ complex matrices,
have been studied by several people \cite{Dubois-Violette, Connes-Lott},
using techniques of non commutative differential geometry. Such
constructions always involve some $\ZZ$-graded differential algebra,
generalizing the usual differential forms. Here, we shall use a differential
calculus that has covariance properties with respect to a finite dimensional
quantum group. Indeed, the algebra of $N \times N$ complex matrices is the
same as the one of a reduced quantum plane with $q^N = 1$ \cite{Weyl}.
Hence we naturally have on ${\mathcal M}$ the action of a quantum group
$\mathcal H$ \cite{Alekseev, Gluschenkov, Coquereaux} and the coaction of
its dual quantum group ${\mathcal F}$. We construct a differential calculus
covariant under the action and coaction of the above mentioned quantum
groups (it is a quotient of the Wess-Zumino complex \cite{Wess-Zumino}) and
study its properties and representation theory (see also \cite{CoGaTr}).
We consider the star operations in ${\mathcal H}$, ${\mathcal F}$ and the
covariant one in ${\mathcal M}$. Working with star representations, the
above mentioned star structures lead to an invariant scalar product in
${\mathcal M}$.

Although our ultimate interest is to build a gauge theory on $S$ with some
invariance property with respect to a quantum group action, the construction
of a Lagrangian will not be considered in the present work. In any case it
is interesting to remark that a general one form on $S$ involves a vector
field $a_\mu$ and two scalar fields $\phi^x$ and $\phi^y$, all of them
valued in ${\mathcal M}\: (=M_N(\CC))$.


\section{The space of complex matrices as a reduced quantum plane}
\label{sec:red-q-plane}

It has been known for a long time \cite{Weyl} that the algebra of
$N \times N$ matrices can be generated by two elements $x$ and $y$ with the
relations
\begin{equation}
xy = q yx
\label{2.1}
\end{equation}

\begin{equation}
x^N = y^N = \one \ ,
\label{2.2}
\end{equation}
where $q$ denotes a $N$-th root of unity ($q \neq 1$) and $\one$ is
the unit matrix. From now on we consider the case $N = 3$. For this case
we have the following explicit representation of $x$ and $y$:
$$
x = \pmatrix{1 & 0 & 0 \cr 0 & q^{-1} & 0 \cr 0 & 0 & q^{-2}} \ ,
\qquad
y = \pmatrix{0 & 1 & 0 \cr 0 & 0 & 1 \cr1 & 0 & 0} \ .
$$
They generate the algebra of $3 \times 3$ complex matrices. The algebra
generated by abstract elements $x,y$ with the relations (\ref{2.1}) is
called the quantum plane ${\CC_q}$, and adding relations (\ref{2.2}) leads
to the reduced quantum plane ${\mathcal M} \doteq M_3(\CC)$. This last
algebra has the following basis as a vector space of dimension nine:
$\{x^r y^s \; ; \;\; r,s=0,1,2\}$.


\section{Quantum group coaction on $\mathcal M$}
\label{sec:q-group-coaction}

Consider the following transformations between ``coordinate functions'',
\begin{equation}
\delta_L \pmatrix{x \cr y} = \pmatrix{a & b \cr c & d}
                             \otimes \pmatrix{x \cr y} \doteq
         \pmatrix{x' \cr y'} \qquad \text{left coaction} \ ,
\label{3.1}
\end{equation}
and
\begin{equation}
\delta_R \pmatrix{x & y} = \pmatrix{x & y} \otimes
                           \pmatrix{a & b \cr c & d} \doteq
         \pmatrix{\tilde x & \tilde y} \qquad \text{right coaction} \ .
\label{3.2}
\end{equation}
These coactions extend to the whole of $\mathcal M$ using the homomorphism
property
$\delta(fg) = \delta (f) \delta(g) \;\;$ $(f,g \in {\mathcal M})$,
for both $L$ and $R$ coactions. The elements $a,b,c,d$ should
satisfy\footnote{
This method of obtaining the product relations was introduced in
\cite{Ogievetsky}.
}
an algebra such that
\begin{eqnarray}
   \delta_L (xy - q yx) &=& 0 \ ,
\label{3.3} \\
   \delta_R (xy - q yx) &=& 0 \ .
\label{3.4}
\end{eqnarray}
Hence one obtains
\begin{equation}
\begin{array}{ll}
   qba = ab \qquad & qdb = bd \\
   qca = ac        & qdc = cd \\
   cb = bc         & ad-da = (q-q^{-1})bc \ , \\
\end{array}
\end{equation}
which are the product relations of what is called $Fun(GL_q(2))$.
The element ${\mathcal D} \doteq da - q^{-1}bc = ad - qbc $ is a central
element (it commutes with all the elements of $Fun(GL_q(2))$); it is
called the q-determinant and we set it equal to $1$, getting
$Fun(SL_{q}(2,\CC))$.

This algebra is a Hopf algebra with the following structure:

\begin{itemize}
\item Coproduct:

     $\Delta a = a \otimes a + b \otimes c$, \qquad
     $\Delta b = a \otimes b + b \otimes d$, \qquad
     $\Delta c = c \otimes a + d \otimes c$, \qquad
     $\Delta d = c \otimes b + d \otimes d$ \qquad \\
     ($\Delta(AB) = \Delta A \Delta B$)

\item Antipode:

     $Sa = d$, \qquad
     $Sb = - q^{-1} b$, \qquad
     $Sc = - q c$, \qquad
     $Sd = a$ \qquad \\
     ($S(uv) = S(v) S(u)$)

\item Counit:

     $\epsilon(a)=1$, \qquad
     $\epsilon(b)=0$, \qquad
     $\epsilon(c)=0$, \qquad
     $\epsilon(d)=1$.

\end{itemize}

If we further include the relations
$$
x'^3 = \one \ , \qquad {\tilde x}^3 = \one \ , \qquad
   y'^3 = \one \ , \qquad {\tilde y}^3 = \one \ ,
$$
we should impose
$$
a^3 = \one \ , \qquad d^3 = \one \ , \qquad b^3 = c^3 = 0 \ ,
$$
that define an ideal and coideal. The quotiented Hopf algebra
so obtained, $\mathcal F$, has dimension $27$.


\section{Quantum group action on $\mathcal M$}
\label{sec:q-group-action}

The dual of the Hopf algebra $\mathcal F$ can be obtained as a quotient of
$U_q(sl(2))$, the dual of $Fun(SL_{q}(2,\CC))$. In order to fix notation
and conventions we recall the defining relations of $U_q(sl(2))$ in terms of
its generators $X_+ , X_- $ and $K$.

\begin{itemize}
\item Multiplication:

     $K X_{\pm} = q^{\pm 2} X_{\pm} K$ , \qquad
     $\left[ X_+, X_- \right] = {K - K^{-1} \over q - q^{-1}}$

\item Comultiplication:

     $\Delta X_+ \doteq X_+ \otimes \one + K \otimes X_+ $ , \qquad
     $\Delta X_- \doteq X_- \otimes K^{-1} + \one \otimes X_-$ \\
     $\Delta K   \doteq K \otimes K$ , \qquad
     $\Delta K^{-1} \doteq K^{-1} \otimes K^{-1}$

\item Antipode:

     $S \one = \one$, \qquad
     $S K = K^{-1}$, \qquad
     $S K^{-1} = K$, \qquad
     $S X_+ = - K^{-1} X_+$, \qquad
     $S X_- = - X_- K$.

\item Counit:

     $\epsilon (\one ) = 1$,
     $\epsilon (K) = 1$,
     $\epsilon (K^{-1}) = 1$,
     $\epsilon (X_\pm ) =0$.

\end{itemize}

\noindent The dual $\mathcal H$ of $\mathcal F$ is obtained by taking the
quotient of $U_q(sl(2))$ with the (Hopf) ideal and coideal defined by
$$
   X_+^3 = X_-^3 = 0 \ , \;\; K^3 = \one \ .
$$

The duality relations between $\mathcal F$ and $\mathcal H$ are given
explicitly by the pairing between generators:
\begin{equation}
\begin{array}{llll}
   <K,a> = q    & <K,b> = 0    & <K,c> = 0    & <K,d> = q^{-1} \\
   <X_+ ,a> = 0 & <X_+ ,b> = 1 & <X_+ ,c> = 0 & <X_+ ,d> = 0   \\
   <X_- ,a> = 0 & <X_- ,b> = 0 & <X_- ,c> = 1 & <X_+ ,d> = 0
\end{array}
\end{equation}
This pairing interchanges multiplication and comultiplication via the
relations
\begin{eqnarray}
   < X_1 X_2, u > = < X_1 \otimes X_2, \Delta u > , \qquad
   <\Delta X, u_1 \otimes u_2 > = < X, u_1 u_2 > \ .
\end{eqnarray}
Since $\mathcal F$ coacts on $\mathcal M$ there is a natural definition of
an action of $\mathcal H$ on $\mathcal M$. Using the pairing and the
{\em right} coaction $\delta_R$, the {\em left} action of $\mathcal H$ on
$\mathcal M$ is defined by
\begin{equation}
h(z) \equiv h.z \doteq (\id \otimes <h,\cdot >) \circ \delta_R z \ ;
     \qquad h \in {\mathcal H}, z \in {\mathcal M} \ .
\label{4.1}
\end{equation}
This implies, in particular, that
\begin{eqnarray}
   h (\one) & \doteq & \epsilon(h) \one \nonumber \\
   h (zw) & \doteq & \Delta(h) (z \otimes w)
\end{eqnarray}
Several properties of this action (with other conventions) have been studied
in \cite{Dabrowski}. Explicitly, one gets
\begin{equation}
\begin{tabular}{l||ccc}
          & $K$        & $X_+$     & $X_-$      \\
\hline
\hline
$x^2$     & $q^2 x^2 $ & $0$       & $-q^2 xy$  \\
$xy$      & $xy$       & $qx^2 $   & $q y^2$    \\
$y^2$     & $q y^2$    & $-q^2 xy$ & $0$        \\
\hline
$x$       & $q x$      & $0$       & $y$        \\
$y$       & $q^2 y$    & $x$       & $0$        \\
$x^2 y^2$ & $x^2 y^2$  & $-qy$     & $-qx$      \\
\hline
$x^2 y$   & $q x^2 y$  & $q^2\one$ & $-xy^2$    \\
$xy^2$    & $q^2 xy^2$ & $-x^2 y$  & $q^2 \one$ \\
$\one$    & $\one $    & $0$       & $0$        \\
\end{tabular}
\label{4.2}
\end{equation}

We conclude this section recalling some facts about the representation
theory of $\mathcal H$.

The left regular representation of $\mathcal H$ was studied in
\cite{Alekseev}. It is the same \cite{Coquereaux} as the $27$-dimensional
non semisimple Hopf algebra{\footnote{
In the next name, $\Lambda^2$ is the Grassmann algebra with two generators
and ${}_0$ denotes the even part of the graded matrix algebra
$M_{2\vert 1}$.}
}
$M_3 \oplus (M_{2\vert 1}(\Lambda^2))_0$ whose elements are given
explicitly by
\begin{equation}
{\mathcal H} = m_3 \oplus
   \pmatrix{ \alpha_{11} + \beta_{11} \theta_1 \theta_2 &
      \alpha_{12} + \beta_{12} \theta_1 \theta_2 &
      \gamma_{13} \theta_1 + \delta_{13} \theta_2 \cr
      \alpha_{21} + \beta_{21} \theta_1 \theta_2 &
      \alpha_{22} + \beta_{22} \theta_1 \theta_2 &
      \gamma_{23} \theta_1 + \delta_{23} \theta_2 \cr
      \gamma_{31} \theta_1 + \delta_{31} \theta_2 &
      \gamma_{32} \theta_1 + \delta_{32} \theta_2 &
      \alpha_{33} + \beta_{33} \theta_1 \theta_2 } \ ,
\label{4.3}
\end{equation}
where the $\alpha$'s, $\beta$'s and the coefficients of the $3 \times 3$
matrix $m_3$ are complex numbers. $\theta_1$ and $\theta_2$ denote Grassmann
variables.

The representations of $\mathcal H$ can be obtained from the action of the
left regular representation matrix of a generic element on its own
columns. The first three columns of the $6 \times 6$ matrix (\ref{4.3})
give equivalent representations that correspond to a three dimensional
irreducible representation denoted by $3_i$. The first two columns of the
$(M_{2\vert 1}(\Lambda^2))_0$ give equivalent representations of dimension
six denoted by $6_e$, while the last column leads to another six-dimensional
representation denoted by $6_o$. These two $6$-dimensional inequivalent
representations are indecomposable but not irreducible, their lattice of
subrepresentations \cite{Coquereaux} is
$$
\begin{diagram}
  0 & \rTo & 2 & \rTo & \{ 3_e^\lambda \} & \rTo & 4_e & \rTo & 6_e \\
  0 & \rTo & 1 & \rTo & \{ 3_o^\lambda \} & \rTo & 5_o & \rTo & 6_o \\
\end{diagram}
$$

\noindent where (column vectors):
$$
\begin{tabular}{ll}
   $ 6_e = (\alpha + \beta \theta_1 \theta_2,
      \alpha' + \beta' \theta_1 \theta_2,
      \gamma \theta_1 + \delta \theta_2) $                            &
   $ 6_o = (\gamma \theta_1 + \delta \theta_2,
      \gamma' \theta_1 + \delta' \theta_2,
      \alpha + \beta \theta_1 \theta_2) $                             \\
   $ 4_e = (\beta \theta_1 \theta_2, \beta' \theta_1 \theta_2,
      \gamma \theta_1 + \delta \theta_2) $                            &
   $ 5_o = (\gamma \theta_1 + \delta \theta_2,
      \gamma' \theta_1 + \delta' \theta_2, \beta \theta_1 \theta_2) $ \\
   $ 3_e^\lambda = (\beta \theta_1 \theta_2, \beta' \theta_1 \theta_2,
      \gamma \theta_\lambda) $                                        &
   $ 3_o^\lambda = (\gamma \theta_\lambda, \gamma' \theta_\lambda,
      \beta \theta_1 \theta_2) $                                      \\
   $ 2 = (\beta \theta_1 \theta_2, \beta' \theta_1 \theta_2, 0) $     &
   $ 1 = ( 0, 0, \beta \theta_1 \theta_2) $                           \\
\end{tabular}
$$

\noindent Here $\alpha, \beta, \gamma, \delta \in \CC$,
$\theta_\lambda = \lambda_1 \theta_1 + \lambda_2 \theta_2$, and
$\lambda_1, \lambda_2 \in \CC$. The families of representations
$\{ 3^\lambda \}$ are actually parametrized by
$\lambda \doteq {\lambda_1 \over \lambda_2} \in CP^1$.


\section{Stars on $\mathcal F$, $\mathcal H$ and $\mathcal M$}
\label{sec:stars}

For the case $|q| = 1 \ (q \neq 1)$ there is only one $*$-Hopf structure on
$Fun(SL_{q}(2,\CC))$ \cite{Chari-Pressley}, up to equivalences{\footnote
{
Two star structures, $\dagger$ and $*$, over a (Hopf) algebra $H$ are
equivalent if there exists a (Hopf) automorphism $\alpha$ such that
$\alpha(a^{\dagger}) = (\alpha(a))^* \: , \ \forall a \in H$. The Hopf
automorphisms of $\mathcal F$ and $\mathcal H$ are \cite{Chari-Pressley}:

\begin{eqnarray*}
{\mathcal F}: &\qquad & a \to a ,\;\; d \to d ,\;\;
                       b \to \alpha b ,\;\; c \to \alpha^{-1} c , \qquad
                       \text{where} \alpha \in \CC \\
{\mathcal H}: &\qquad & K \to K ,\;\; X_+ \to \beta X_+ ,\;\;
                       X_- \to \beta^{-1} X_- , \qquad
                       \text{where} \beta \in \CC
\end{eqnarray*}

\noindent Using the above notion of equivalence of star structures we get
the following family of equivalent stars:
\begin{eqnarray*}
{\mathcal F}: &\qquad & a^{*_u} = a ,\;\; d^{*_u} = d ,\;\;
                       b^{*_u} = u b ,\;\; c^{*_u} = u^{-1} c , \qquad
                       \text{with} |u|=1 \\
{\mathcal H}: &\qquad & K^{*_u} = K ,\;\; X_+^{*_u} = -q^{-1} u X_+ ,
                       \;\; X_-^{*_u} = -q u^{-1} X_- \ .
\end{eqnarray*}
}
}. It is given by
\begin{equation}
a^* = a \ ,\;\; b^* = b \ ,\;\; c^* = c \ ,\;\; d^* = d \ .
\end{equation}

\noindent We denote this last $*$-Hopf algebra as $Fun(SL_{q}(2,\RR))$, in
correspondence with the classical case. This star is compatible with the
relations $a^3 = d^3 = \one, \ b^3 = c^3 = 0$, and therefore extends to the
quotient of $Fun(SL_{q}(2,\RR))$ that we continue to denote by $\mathcal F$.

Using the pairing between $\mathcal F$ and $\mathcal H$ we can obtain the
corresponding star in $\mathcal H$. In order to get a $*$-Hopf structure in
the dual $\mathcal H$ of $\mathcal F$ one should have \cite{Chari-Pressley} 
\begin{equation}
< h^* , a > = \overline{<h , (Sa)^* >} \ ,
\label{5.00}
\end{equation}
where the bar means complex conjugation. This leads to the following star
in $\mathcal H$,
\begin{equation}
{X_+}^* = - q^{-1} X_+ \ , \;\; {X_-}^* = - q X_- \ , \;\;K^* = K \ .
\end{equation}

Next we consider the determination of a star in $\mathcal M$. The basic
requirement is that the coaction $\delta \doteq \delta_R$ should be a
$*$-homomorphism, \ie
\begin{equation}
(\delta z)^* = \delta(z^*) \ , \;\; \text{for any} z \in {\mathcal M} ,
	     \;\;\; \text{where} (A \otimes B)^* = A^* \otimes B^* \ .
\label{5.0}
\end{equation}
Pairing the $\mathcal F$ component of this equation with some 
$h \in {\mathcal H}$, and recalling (\ref{4.1}) and (\ref{5.00}), we see 
that the dual of condition (\ref{5.0}) is
\begin{equation}
h(z^*) = \left[ (S h )^* x \right]^* \ .
\label{5.1}
\end{equation}
It is not difficult to verify that the only family of stars compatible with
(\ref{5.1}) is
\begin{equation}
x^* = \alpha x \ ,\;\; y^* = \alpha y \ , \qquad \text{with} |\alpha|=1 \ ;
\end{equation}
moreover, we chose{\footnote
{
Changing the star in $\mathcal F$ and $\mathcal H$ to an equivalent one as
in the previous footnote, corresponds to a change in the family of stars in
$\mathcal M$ given by
\begin{equation}
x^{*_u} = \alpha \, x \ ,\;\; y^{*_u} = u \alpha \, y \ , \qquad
          |\alpha|=1.
\end{equation}
}
} $\alpha = 1$.


\section{Quantum groups and invariant scalar products}
\label{sec:scalar-products}

It is not straightforward to define the notion of an invariant scalar
product on a representation space $V$ of a quantum group $\mathcal H$. The
basic idea is that taking scalar products should in some way commute with
the action of the quantum group. This amounts to say that
\begin{equation}
   h \circ (\cdot,\cdot) = (\cdot,\cdot) \circ {\tilde h}
\label{6.1}
\end{equation}
Here, $(\cdot , \cdot )$ denotes the scalar product on $V$ (antilinear in
the first variable and linear in the second).  The action of $\mathcal H$
on $\CC$ is the trivial one, given by $h(\alpha) = \epsilon(h) \alpha$,
so that
\begin{equation}
   (h \circ (\cdot,\cdot))(u \otimes v) = h((u,v)) = \epsilon(h) (u,v) \ .
\label{6.2}
\end{equation}
Finally, ${\tilde h}$ is an action of $\mathcal H$ on $V \otimes V$ to be
determined.

The natural action of $\mathcal H$ on $V \otimes V$ is given by
\begin{equation}
   h.(u \otimes v) = \Delta h (u \otimes v) = (h_1 u) \otimes (h_2 v) \ ,
\end{equation}
however this action can not be the one involved in the r.h.s. of requirement
(\ref{6.1}). This is so since the l.h.s. of this equation ---as given by
(\ref{6.2})--- has to be linear under the replacement $h \to \alpha h$,
with $\alpha \in \CC$, but the right hand side would be linear if we attach
$\alpha$ to $h_2$ and antilinear if we attach $\alpha$ to $h_1$. If there
is a star structure in $\mathcal H$ one could try to solve this problem by
defining ${\tilde h}(u \otimes v)$ as $h_1^* u \otimes h_2 v$. This would
fix the sesquilinearity problem. Unfortunately, this choice for ${\tilde h}$
would not have a definite homomorphism behaviour for the product in
$\mathcal H$ whereas the $\epsilon$ on the r.h.s. of (\ref{6.2}) (and thus
on the l.h.s. of (\ref{6.1})), is a homomorphism. The solution to this last
problem is to employ the natural antihomomorphism given by the antipode.
However there are two options,

\begin{equation}
   {\tilde h^1} (u \otimes v) = (S h_1 )^* u \otimes h_2 v
\end{equation}
or
\begin{equation}
   {\tilde h^2} (u \otimes v) = S(h_1^*) u \otimes h_2 v \ .
\end{equation}
Replacing these for $\tilde h$ in (\ref{6.1}), we get the corresponding
conditions of invariance, which are
\begin{equation}
   \epsilon(h) (u,v) = ((S h_1 )^* u , h_2 v)
\label{6.3}
\end{equation}
and
\begin{equation}
   \epsilon(h) (u,v)=  (S(h_1^*) u , h_2 v) \ .
\label{6.4}
\end{equation}
We will choose to work with $*$-representations, \ie representations
of $\mathcal H$ such that
\begin{equation}
   (hu, v) = (u, h^* v) \ .
\label{6.5}
\end{equation}
In this case the requirement (\ref{6.3}) is automatically fulfilled since
\begin{equation}
   ((S h_1 )^* u, h_2 v) = (u, S(h_1) h_2 v) =
      (u, m[(S \otimes \id) \Delta h ].v) = \epsilon (h) (u,v) \ .
\end{equation}

\noindent This result contrasts with the later condition (\ref{6.4}) which
is not a consequence of the star representation condition (\ref{6.5}).
Condition (\ref{6.4}) is actually not fulfilled in our case
($V = {\mathcal M}$) unless we choose a vanishing scalar product.

Having decided to work with star representations of $\mathcal H$, the scalar
product will then be automatically invariant in the sense of requirement
(\ref{6.3}), \ie for ${\tilde h}\doteq {\tilde h^1}$, and non-vanishing.


\subsection{Scalar product for the representation of $\mathcal H$ on
            $\mathcal M$}

As $\mathcal M$ acts on itself by left multiplication, we will also require
that this representation should also be a star representation, \ie
$(z,z') = (\one, z^* z') = (z'^* z, \one)$ for any $z,z' \in {\mathcal M}$.
This last requirement implies that it is enough to know the scalar products
of the form $(\one,z), \;\; z \in {\mathcal M}$. The condition of star
representation for the action of $\mathcal H$ on these particular scalar
products is
\begin{equation}
   (h^* .\one , z) = (\one, h.z) \ , \qquad h \in {\mathcal H}
\label{6.1.1}
\end{equation}
Taking $z = x^r y^s$ and $h = K$ one concludes from (\ref{4.2}) and
(\ref{6.1.1}) that the scalar products of the form $(\one, x^r y^s)$
vanish if $r \neq s$. Taking $h = X_+$ one gets that the only one that
does not vanish is
\begin{equation}
   (\one, x^2 y^2) \neq 0 \ .
\label{6.1.2}
\end{equation}
Taking $h = X_-$ gives no further condition. Recalling that we work in a
star representation for the action of $\mathcal M$, we get from
(\ref{6.1.2}) eight additional nonvanishing scalar products. We fix their
value by choosing
\begin{equation}
   (xy,xy) = 1 \ .
\end{equation}


\section{Covariant differential calculus on $\mathcal M$}
\label{sec:diff-calculus}

We will use the covariant differential calculus of \cite{Wess-Zumino}.
This calculus is built taking the Manin dual of the quantum plane
\cite{Manin} as the algebra of differentials, the algebra between the
differentials and quantum plane variables been obtained requiring
quadraticity, covariance and consistency. These relations are

\begin{equation}
\begin{tabular}{ll}
$xy = q\, yx$ \ ,             &                                   \\
$x\, dx = q^2 dx\, x$ \ ,     &
          $x\, dy =q\, dy\, x +(q^2-1)\, dx\, y$ \ ,              \\
   $y\, dx = q\, dx\, y $ \ , & $y\, dy = q^2 dy\, y$             \\
$dx^2 = dy^2 = 0$ \ ,         & $dx \, dy + q^2 dy \, dx = 0$ \ . \\
\end{tabular}
\label{7.1}
\end{equation}

\noindent In order to extend this algebra to a differential calculus on
the reduced quantum plane $\mathcal M$, we should check that the ideal
employed in defining $\mathcal M$ is a differential one. Indeed,
it is simple to verify that
\begin{equation}
   d (x^3) = 0 = d (y^3) \ .
\label{7.2}
\end{equation}
This differential algebra, $\Omega_{WZ}({\mathcal M})$, is graded with
the decomposition
\begin{equation}
   \Omega_{WZ}({\mathcal M}) =
      \displaystyle{\bigoplus_{n=0}^2} \ \Omega^n_{WZ}({\mathcal M}) \ ,
\end{equation}
where
\begin{equation}
\Omega^0_{WZ}({\mathcal M}) = {\mathcal M} \ , \qquad
\Omega^1_{WZ}({\mathcal M}) = {\mathcal M} \: dx \oplus
                              {\mathcal M} \: dy \ , \qquad
\Omega^2_{WZ}({\mathcal M}) = {\mathcal M} \: dx \, dy \ .
\end{equation}


\section{The action of $\mathcal H$ on $\Omega_{WZ}({\mathcal M})$}
\label{sec:action-on-omega}


\subsection{The action of $\mathcal H$ on
        $\Omega^0_{WZ}({\mathcal M}) = {\mathcal M}$}

This action is given in (\ref{4.2}). From that table and the representation
theory of $\mathcal H$ (see the summary at the end of
Section~\ref{sec:q-group-action}), we see that $x^2, xy, y^2$ span the $3_i$
representation. Moreover, $x,y,x^2 y^2$ span one of the $3$-dimensional
representations $3_e^\lambda \subset 6_e$ (which we will call simply $3_e$).
Finally, $x^2 y, \one, x y^2$ span one of the $3$-dimensional
representations $3_o^\lambda \subset 6_o$ (which we will call simply $3_o$).
The above mentioned action is represented schematically in the following
figures, where dashed arrows stand for the action of $X_-$ and continuous
ones for $X_+$:

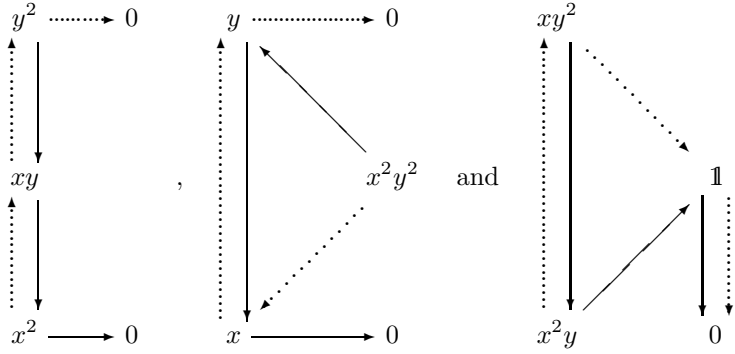
\vspace{0.3cm}
$$
\begin{diagram}
  y^2           & \rDotsto 0 \\
  \uDotsto \dTo &            \\
  xy            &            \\
  \uDotsto \dTo &            \\
  x^2           & \rTo 0
\end{diagram}
\hspace{0.5cm} \hbox{,} \hspace{0.5cm}
\begin{diagram}
  y             & \rDotsto  & 0       \\
                & \luTo     &         \\
  \uDotsto \dTo &           & x^2 y^2 \\
                & \ldDotsto &         \\
  x             & \rTo      & 0
\end{diagram}
\hspace{0.5cm} \hbox{and} \hspace{0.5cm}
\begin{diagram}
  xy^2          &           &               \\
                & \rdDotsto &               \\
  \uDotsto \dTo &           & \one          \\
                & \ruTo     & \dTo \dDotsto \\
  x^2 y         &           & 0
\end{diagram}
$$
\vspace{0.3cm}

\noindent Note that $3_e$ contains an irreducible $2$-dimensional
representation spanned by $x,y$, and that $3_o$ contains an irreducible
$1$-dimensional representation spanned by $\one$.


\subsection{The action of $\mathcal H$ on $\Omega^1_{WZ}({\mathcal M})$}

First we consider the action of $\mathcal H$ on the Manin dual
$\mathcal M^!$ of $\mathcal M$, spanned by $dx$ and $dy$. We note that the
right coaction of $\mathcal F$ on $\mathcal M^!$ is the same as in
(\ref{3.2}), replacing $x,y$ by $dx, dy$ \cite{Manin}. Indeed, one can
alternatively define $Fun(SL_q(2,\CC))$ by requiring only (\ref{3.4}) and
the invariance of the algebraic relations of the Manin dual $\mathcal M^!$
under the above mentioned coaction. The left action of $\mathcal H$ on
$\mathcal M^!$ is again defined by the relation (\ref{4.1}), hence we get:

$$
   K dx = q dx \:, \ K dy = q^{-1} dy \ ; \;\;\;
   X_+ dx = 0  \:, \ X_+ dy = dx \ ; \;\;\;
   X_- dx = dy \:, \ X_- dy = 0 \ .
$$

\noindent This action corresponds to the irreducible representation $2$ of
$\mathcal H$. Therefore, in order to know the transformation properties of a
1-form we should decompose the tensor products $3_i \otimes 2$,
$3_e \otimes 2$ and $3_o \otimes 2$ into direct sums of the indecomposable
representations. These tensor products are constructed using the
coproduct. The result for these decompositions follows.

\begin{itemize}
\item
   The case $3_o \otimes 2 = 3_i \oplus 3_e$.

\vspace{0.3cm}
\hspace{0.8cm}
\begin{tabular}{l||ccc}
              & $K$               & $X_+$          & $X_-$ \\
\hline
\hline
$x^2 y \, dx$ & $q^2 x^2 y \, dx$ & $q^2\, dx$     &
                $-q^2 x y^2 \, dx + x^2 y \, dy$           \\
$xy^2 \, dx $ & $x y^2 \, dx$     & $-x^2 y \, dx$ &
                $q dx + xy^2 \, dy$                        \\
$\one \, dx$  & $q \, dx$         & $0$            & $dy$  \\
\hline
$x^2 y \, dy$ & $x^2 y \, dy$     & $q^2 \, dy + q x^2 y \, dx$      &
                $-q x y^2 \, dy$                           \\
$x y^2 \, dy$ & $q x y^2 \, dy$   & $-x^2 y \, dy + q^2 x y^2 \, dx$ &
                $dy$                                       \\
$\one \, dy$  & $q^2 \, dy$       & $dx$           & $0$   \\
\end{tabular}
\vspace{0.3cm}

This indeed gives $3_i \oplus 3_e$, since, up to multiplicative
factors,

\vspace{0.3cm}
$$
\begin{diagram}
  q \, x^2 y \, dx - dy          & \rTo 0     \\
  \uTo^{} \dDotsto^{}            &            \\
  -x^2 y \, dy + q^2 x y^2 \, dx &            \\
  \uTo^{} \dDotsto^{}            &            \\
  -q \, dx + q \, x y^2 \, dy    & \rDotsto 0 \\
\end{diagram}
\hspace{1.0cm} \oplus \hspace{1.0cm}
\begin{diagram}
  dy                  & \rDotsto  & 0                              \\
                      & \luTo     &                                \\
  \uDotsto^{} \dTo^{} &           & x^2 y \, dy + q \, x y^2 \, dx \\
                      & \ldDotsto &                                \\
  dx                  & \rTo      & 0
\end{diagram}
$$
\vspace{0.3cm}

Notice that $3_e$ contains an irreducible $2$-dimensional representation
spanned by $\{ dx, dy \}$.

\item
   The case $3_e \otimes 2 = 3_i \oplus 3_o$.

\vspace{0.3cm}
\hspace{0.8cm}
\begin{tabular}{l||ccc}
                & $K$               & $X_+$        & $X_-$     \\
\hline
\hline
$x \, dx$       & $q^2 x \, dx$     & $0$          &
                  $q^2 y \, dx + x \, dy$                      \\
$y \, dx$       & $y \, dx$         & $x \, dx$    & $ y\, dy$ \\
$x^2 y^2 \, dx$ & $q x^2 y^2 \, dx$ & $-q y \, dx$ &
                  $-x \, dx + x^2 y^2 \, dy$                   \\
\hline
$x \, dy$      & $x \, dy$           & $q x \, dx$                  &
                 $q y \, dy$                                         \\
$y \, dy$      & $qy \, dy$          & $x \, dy + q^2 y \, dx$      &
                 $0$                                                 \\
$x^2y^2 \, dy$ & $q^2 x^2 y^2 \, dy$ & $-q y \, dy + x^2 y^2 \, dx$ &
                 $-q^2 x \, dy$                                      \\
\end{tabular}
\vspace{0.3cm}

This gives $3_i \oplus 3_o$, since, up to multiplicative
factors,

\vspace{0.3cm}
$$
\begin{diagram}
  y \, dy                  & \rDotsto 0 \\
  \uDotsto^{} \dTo^{}      & \\
  q \, x \, dy + y \, dx   & \\
  \uDotsto^{}  \dTo^{}     & \\
  x \, dx                  & \rTo 0 \\
\end{diagram}
\hspace{1.0cm} \oplus \hspace{1.0cm}
\begin{diagram}
  x^2y^2\, dy - x \, dx &           & \\
                        & \rdDotsto & \\
  \uDotsto^{} \dTo^{}   &           & x\, dy - q \, y \, dx \\
                        & \ruTo     & \dTo \dDotsto \\
  y \, dy - q^2 \, x^2 y^2\, dx &   & 0
\end{diagram}
$$
\vspace{0.3cm}

Note that $3_o$ contains an irreducible $1$-dimensional representation
spanned by $x \, dy - q \, y \, dx$.

\item
   The case $3_i \otimes 2 = 6_e$.

\vspace{0.3cm}
\hspace{0.8cm}
\begin{tabular}{l||ccc}
             & $K$            & $X_+$          & $X_-$                    \\
\hline
\hline
$x^2 \, dx$  & $x^2 \, dx$    & $0$            & $-qxy \, dx + x^2 \, dy$ \\
$xy \, dx$   & $qxy \, dx$    & $qx^2 \, dx$   & $ y^2 \, dx + xy \, dy$  \\
$y^2 \, dx$ & $q^2 y^2 \, dx$ & $-q^2xy \, dx$ & $y^2 \, dy$              \\
\hline
$x^2 \, dy$  & $q x^2 \, dy $ & $q^2 x^2 \, dx$         & $-xy \, dy$     \\
$xy \, dy$   & $q^2 xy \, dy$ & $qx^2 \, dy + xy \, dx$ & $q^2 y^2 \, dy$ \\
$y^2 \, dy$  & $y^2 \, dy$    & $-q^2xy \, dy+q y^2 \, dx$ & $0$          \\
\end{tabular}
\vspace{0.3cm}

This is the six-dimensional indecomposable representation $6_e$ (which is
projective, \cf \cite{Coquereaux}). It contains the four-dimensional
indecomposable representation $4_e$, a family of indecomposables of the
type $3_e^\lambda$, and one irreducible of dimension two, spanned by
$\{ -q^2 x^2 \, dy + xy \, dx, -q xy \, dy + y^2 \, dx\}$.

\vspace{0.3cm}
$$
\begin{diagram}
  0                         &       & x^2\, dx  & \rTo     & 0            \\
  \uTo                      & \ldDotsto &       & \luTo    &              \\
  -q^2 x^2 \, dy + xy \, dx &       &           &          &
                              xy \, dx + x^2 \, dy                        \\
  \uTo \dDotsto             & & \ldDotsto\luTo  &          & \uTo\dDotsto \\
  -q xy \, dy + y^2 \, dx   &       &           &          &
                              xy \, dy + y^2 \, dx                        \\
  \dDotsto                  & \luTo &          & \ldDotsto &              \\
  0                         &       & y^2 \, dy & \rDotsto & 0            \\
\end{diagram}
$$
\vspace{0.3cm}
\end{itemize}


\subsection{The action of $\mathcal H$ on $\Omega^2_{WZ}({\mathcal M})$}

Noting that $\Omega^2_{WZ}({\mathcal M}) = {\mathcal M} \: dx \, dy$
(and thus it is isomorphic to $\mathcal M$), and that $dx \, dy$ is
invariant under the action of $\mathcal H$ we conclude that
$\Omega^2_{WZ}({\mathcal M})$ decomposes in exactly the same
representations as ${\mathcal M}$.


\section{Star structure on $\Omega_{WZ}({\mathcal M})$}
\label{sec:star-on-omega}

We now consider the possible stars for $\Omega_{WZ}({\mathcal M})$. Since we
have an action of $\mathcal H$ on $\Omega_{WZ}({\mathcal M})$, we will look
for stars compatible with this action, \ie the ones satisfying the analog
of (\ref{5.1}) for this case,
\begin{equation}
   h(\omega^*) = ((S h)^* \omega)^* \ , \qquad
        \omega \in \Omega_{WZ}({\mathcal M}) \ .
\label{9.1}
\end{equation}

This determines the action of $\mathcal H$ on
$\Omega^*_{WZ}({\mathcal M})$ ($\simeq \Omega_{WZ}({\mathcal M})$).
Moreover, we will also require the stars to respect the grading of
$\Omega_{WZ}({\mathcal M})$. We already know the star on
$\Omega^0_{WZ}({\mathcal M}) = {\mathcal M}$. For the case of
$\Omega^1_{WZ}({\mathcal M})$ we proceed by writing the following most
general expressions for $(dx)^*$ and $(dy)^*$,
\begin{eqnarray}
   dx^* &=& \alpha^1_{rs} \, x^r y^s dx +
            \alpha^2_{rs} \, x^r y^s dy          \nonumber \\
   dy^* &=& \beta^1_{rs} \, x^r y^s dx + \beta^2_{rs} \, x^r y^s dy \ .
\end{eqnarray}
Replacing this in (\ref{9.1}) we obtain equations for the coefficients
$\alpha$ and $\beta$. Furthermore, imposing that
\begin{equation}
   x^* dx^* = q dx^* x^* \ ,
\end{equation}
which is the star of one of the relations in (\ref{7.1}), it can be seen
that the only solution is
\begin{equation}
   dx^* = dx \ , \qquad dy^* = dy \ .
\end{equation}

We find useful to remark that in spite of being on a non commutative space
this star on $\Omega_{WZ}({\mathcal M})$ shares some properties with
the one for real manifolds. It can be easily checked that for all elements
$z$ in ${\mathcal M}$, one has $(dz)^* = d(z^*)$. Therefore
(with $\phi = z_i \, dx^i \in \Omega_{WZ}^1$), $d(\phi^*) = - (d\phi)^*$.
More generally, $d(\omega^*) = (-1)^p (d\omega)^*$, when
$\omega \in \Omega_{WZ}^p$.

Having a star we can write the most general hermitian one-form on
$\mathcal M$. It is given by a (real) linear combination of the following
hermitian forms:

\begin{eqnarray*}
   \omega_{3i}  &=& \alpha_1 (q x^2 y \, dx - dy) +
      \alpha_2 (q^2 x y^2 \, dx - x^2 y \, dy) +
      \alpha_3 (q dx - q x y^2 \, dy)                           \\
   \omega'_{3i} &=& \alpha'_1 q^2 y\, dy +
      \alpha'_2 (y\, dx + q x\, dy) + \alpha'_3 q^2 x\, dx      \\
   \omega_{3e}  &=& \beta_1 dy + \beta_2 (q x y^2 \, dx + x^2 y \, dy) +
      \beta_3 dx                                                \\
   \omega_{3o}  &=& \beta'_1 q (x\, dx - x^2y^2\, dy) +
      \beta'_2 (q^2 y\, dx - q x\, dy) +
      \beta'_3 q^2 (q^2 x^2 y^2\, dx - y\, dy)                  \\
   \omega_{6e}  &=& \gamma_1 (xy\, dx - q^2 x^2\, dy) +
      \gamma_2 (y^2\, dx - q xy\, dy) +
      \gamma_3 q x^2\, dx + \gamma_4 q y^2\, dy +               \\
   & &\gamma_5 q^2 (xy\, dx + x^2\, dy) +
      \gamma_6 q (xy\, dy + y^2\, dx)                           \\
\end{eqnarray*}

\noindent The coefficients $\alpha_i, \beta_i, \gamma_i$ are arbitrary
real numbers, and the subscripts of the $\omega$'s refer to the
indecomposable representations previously discussed.


\section{Incorporation of space-time}
\label{sec:space-time}

Let $\Lambda$ be the algebra of usual differential forms over a space-time
manifold $M$ (the De Rham complex) and
$\Omega_{WS} \doteq \Omega_{WS}({\mathcal M})$
the differential algebra over the reduced quantum plane introduced in
Section~\ref{sec:diff-calculus}. Remember that
$\Omega_{WZ}^0 = {\mathcal M}$,
$\Omega_{WZ}^1 = {\mathcal M} \: dx + {\mathcal M} \: dy$, and that
$\Omega_{WZ}^2 = {\mathcal M} \: dx \, dy$.
We call $\Xi$ the graded tensor product of these two differential algebras:
$$
   \Xi \doteq \Lambda \otimes \Omega_{WZ}
$$

\begin{itemize}
\item
   A generic element of $\Xi^0 = \Lambda^0 \otimes \Omega_{WZ}^0$ is a
   $3 \times 3$ matrix with elements in $C^\infty(M)$. It can be thought
   as a scalar field valued in $M_3(\CC)$.

\item
   A generic element of
   $\Xi^1 = \Lambda^0 \otimes \Omega_{WZ}^1 \oplus
            \Lambda^1 \otimes \Omega_{WZ}^0$
   is given by a triplet $\omega = \{ a_{\mu}, \phi^x, \phi^y \}$, where
   $a_{\mu}$ determines a one-form (a vector field) on the manifold $M$
   with values in $M_3(\CC)$ (that we can consider as the Lie algebra of the
   Lie group $GL(3,\CC)$), and where $\phi^x, \phi^y$ are $M_3(\CC)$-valued
   scalar fields. Indeed
   $\phi^x (x_{\mu}) \: dx + \phi^y (x_{\mu}) \: dy
            \in \Lambda^0 \otimes \Omega_{WZ}^1$.

\item
   A generic element of
   $\Xi^2 = \Lambda^0 \otimes \Omega_{WZ}^2 \oplus
            \Lambda^1 \otimes \Omega_{WZ}^1 \oplus
            \Lambda^2 \otimes \Omega_{WZ}^0$
   consists of
   \begin{itemize}
   \item
      a matrix-valued $2$-form $F_{\mu \nu} dx^\mu dx^\nu$ on the
      manifold $M$, \ie an element of $\Lambda^2 \otimes \Omega_{WZ}^0$

   \item
      a matrix-valued scalar field on $M$, \ie an element of
      $\Lambda^0 \otimes \Omega_{WZ}^2$

   \item
      two matrix-valued vector fields on $M$, \ie an element of
      $\Lambda^1 \otimes \Omega_{WZ}^1$
   \end{itemize}
\end{itemize}

The algebra $\Xi$ is endowed with a differential (of square zero, of course,
and obeying the Leibniz rule) defined by
$d \doteq d \otimes \id \pm \id \otimes d$. Here $\pm$ is the
(differential) parity of the first factor of the tensor product upon
which $d$ is applied, and the two $d$'s appearing on the right hand side
are the usual De Rham differential on antisymmetric tensor fields and the
differential of the reduced Wess-Zumino complex, respectively.


\section{Concluding remarks}
\label{sec:remarks}

Fundamental interactions are usually described by gauge theories (abelian
or not) that may have two kinds of invariances: global and local symmetries.
In both cases, the geometrical interpretation is clear and is ultimately
described in terms of group actions (gauged or not).

On the other hand, quantum groups ---either specialized at roots of unity
or not--- have been used many times, during the last decade, in the physics
of integrable models and in conformal theories, but not in space-time
physics. There is no good reason (no known reason), however, to discard
quantum groups from the toolbox used to construct four-dimensional classical
---or quantum--- field theories. Several related attempts have been made to
replace the four-dimensional space-time itself by a quantum space on which
a quantum group would act, or to develop the analog of gauge theories for
quantum groups. Our present goal, in contrast, is to develop a new type of
field theories, where global symmetries are described by a quantum group (we
use the adjective ``global'' since we are not at all trying to gauge
---whatever it means--- the action of the quantum group). This program
requires a new understanding of several basic notions, including reality
structure (hermiticity) and differential calculus; this was discussed here.
Other aspects shall be discussed in a forthcoming review paper \cite{CoGaTr}.

Besides casting some light on representation theory of a specific class
of non semi-simple finite dimensional Hopf algebras, the present work
provides a first example of a particular type of generalized differential
forms, defined on a space-time manifold, and having covariance properties
both with respect to the Lorentz group (or any group acting on the chosen
space-time $M$) and with respect to a finite dimensional quantum group
acting on ``internal indices''. More precisely, if $G$ is a Lie group
acting on $M$, and if ${\mathcal U} $ denotes the enveloping algebra of its
Lie algebra, our construction gives an action of the Hopf algebra
${\mathcal U} \otimes {\mathcal H}$ on the differential algebra $\Xi$.

As stated in the Introduction, we did not undertake in this work, the
construction of any kind of Lagrangian model, but we hope that the present
paper will trigger some interesting ideas in that direction.


\section{Acknowledgements}

We are greatly indebted to Oleg Ogievetsky for his comments and discussions.


\end{document}